\documentclass[german,english
              ,twoside
              ,a4paper
              ,10pt
              ]{amsart}
\usepackage{gt-sim}
\usepackage[T1]{fontenc}
\usepackage[latin1]{inputenc}
\usepackage{babel}
\usepackage{hyperref}
\usepackage{cas-math}
\usepackage{cas-paper}
\usepackage{color}
\usepackage[arrow,cmtip,matrix]{xy}
\usepackage{epsfig}

\definecolor{draftmargin}{rgb}{0.7,0,0}

\newcommand\draftcomment[1]{}

\newenvironment{acknowledgements}{\subsection{Acknowledgements}}{}

\let\phi\varphi
\def\opp#1{{#1^o}}
\def\real{\abs}
\def\ctimes{\times}
\DeclareMathOperator{\ind}{ind}
\DeclareMathOperator{\Ind}{Ind}
\DeclareMathOperator{\Int}{Int}
\DeclareMathOperator{\lk}{lk}

\begin{document}

\title[Small models of graph
  colouring manifolds and $\Hom(C_5,K_n)$]%
{Small models of graph
  colouring manifolds\\and the Stiefel manifolds $\Hom(C_5, K_n)$}
\author{Carsten Schultz}
\address{Institut für Mathematik, MA 6-2\\
Technische Universität Berlin\\
D-10623 Berlin, Germany}
\email{carsten@codimi.de}
\thanks{This research was supported by the
\foreignlanguage{german}{Deutsche Forschungsgemeinschaft} within the
European graduate program ``Combinatorics, Geometry, and Computation''
(No.~GRK~588/2)}

\date{October 2005, introduction updated in July 2006}

\sloppy

\begin{abstract}
We show Péter Csorba's conjecture that the graph homomorphism complex
$\Hom(C_5,K_{n+2})$ is homeomorphic to a Stiefel manifold, the space
of unit tangent vectors to the $n$-dimensional sphere.  For this a
general tool is developed that allows to replace the complexes
$\Hom(G, K_n)$ by smaller complexes that are homeomorphic to them
whenever $G$ is a graph for which those complexes are manifolds.  The
equivariant version of Csorba's conjecture is proved up to homotopy.

We also study certain subdivisions of simplicial manifolds that are
related to the interval poset of their face posets and their
connection with geometric approximations to diagonal maps.
\end{abstract}

\maketitle


\section{Introduction}
\label{sec:intro}
The study of properties of graphs, especially their chromatic number,
through topological spaces associated to the graphs began with Lovász'
proof of Kneser's conjecture~\cite{lovasz}.  Recently the focus has
shifted to the homomorphism complex $\Hom(G,H)$ associated to two
graphs $G$ and~$H$, see~\cite{babson-kozlov-i}.  It is a cell complex
whose vertices are the graph homomorphisms from $G$ to~$H$ and whose
topology captures the way in which homomorphisms can be transformed
into each other by local changes; its cells correspond to
multi-homomorphisms from $G$ to~$H$, functions which assign to every
vertex of $G$ a set of vertices of~$H$ such that every choice
function for it is a homomorphism, see \prettyref{def:Hom}.  In
particular Babson and Kozlov have proved a conjecture of Lovász that
states that if for a graph~$G$ and an $r\ge1$ the complex
$\Hom(C_{2r+1}, G)$ is $(n-1)$-connected then $G$ is not
$(n+2)$-colourable~\cite{babson-kozlov-ii,codd}.  The proof uses the
functoriality of $\Hom$ and topological properties of the complexes
$\Hom(C_{2r+1}, K_{n+2})$.

Despite these advances, the homotopy or homeomorphism types of very
few of the complexes $\Hom(G,H)$ are known, even in the case where $H$
is a complete graph and graph homomorphisms become colourings.  The
cohomology groups of the spaces $\Hom(C_m, K_{n+2})$ have been
calculated in~\cite{kozlov-cycles}.  Susequent to this and the current
work it has been established in \cite{hom-loop} that the colimit of a
diagram of all $\Hom(C_m, K_{n+2})$ with even~$m$ is homotopy
equivalent to the free loop space of~$\Sphere^n$ while the colimit
over all odd~$m$ is homotopy equivalent to the space of all loops
in~$\Sphere^n$ which are equivariant with repsect to the antipodal
actions on $\Sphere^1$ and~$\Sphere^n$.

Among
the spaces $\Hom(C_m, K_{n+2})$, the spaces
$\Hom(C_5, K_{n+2})$ are special, because they are manifolds.  Spaces
$\Hom(G, K_n)$ which are manifolds, graph colouring manifolds, have
been studied by Csorba and Lutz~\cite{csorba-lutz}.  Based on
cohomology and index calculations, Péter Csorba has conjectured that
$\Hom(C_5, K_{n+2})$ is homeomorphic to the Stiefel manifold of
orthonormal $2$-frames in~$\R^{n+1}$~\cite[Conj.~4.8]{csorba-thesis}.
A proof of this is the main result of this article.
\begin{thm*}[\ref{thm:result}]
Let $n\ge0$.  Then there is a homeomorphism
\begin{equation}\label{eq:main}
\Hom(C_5,
K_{n+2})\homeo\set{(x,y)\in\Sphere^n\times\Sphere^n\colon\kron xy=0}.
\end{equation}
\end{thm*}
The case $n=0$, $\Hom(C_5, K_2)=\emptyset$, is trivial, since $C_5$
admits no $2$-colouring. The case $n=1$, $\Hom(C_5,
K_3)\homeo\Sphere^1\times\Sphere^0$, is easily checked.  The cases
$n=2$ and $n=3$ are proved in \cite{csorba-lutz}, $\Hom(C_5,
K_4)\homeo\R P^3$ by direct construction, and $\Hom(C_5,
K_5)\homeo\Sphere^3\times\Sphere^2$ by a one week computer
calculation.

Our proof consists of two distinct parts, and this article is
structured accordingly.  The first part consists of finding a smaller
model for~$\Hom(C_5, K_{n+2})$ and is presented in
\prettyref{sec:res}.  If $I$ is a set of vertices of a graph~$G$, then
one can define a complex~$\Hom_I(G, H)$ that is derived from $\Hom(G,
H)$ by identifying graph homomorphisms that differ only on vertices
in~$I$.  If $I$ is an independent set of vertices, then $\Hom_I(G, H)$
is homotopy equivalent to $\Hom(G, H)$.  This has first been proved
and used in \cite[Sec.~2.8]{csorba-thesis}.  We show that for graph
colouring manifolds it does indeed yield homeomorphic complexes.
\begin{lem*}[\ref{thm:hom-mf}, \ref{lem:mf-homeo}]
Let $G$ be a graph and $I$ an independent set of vertices of~$G$.  If
$\Hom(G, K_n)$ is a manifold for all~$n$, then 
$\Hom(G, K_n)\homeo\Hom_I(G, K_n)$ for all~$n$.
\end{lem*}
We hope that this result will find further applications.  At least the
fact alone that $\Hom_I(G, H)$ is a complex with a smaller number of
vertices than $\Hom(G, H)$ should make the determination of the
homeomorphism types of more graph colouring manifolds accessible to
computer calculations just Péter Csorba has used the homotopy
equivalence of $\Hom(G, H)$ and $\Hom_I(G, H)$ to, among other things,
simplify computer calculations of cohomology groups of homomorphism
complexes.  In the current work, however, this construction is used to
expose the structure of $\Hom(C_5, K_{n+2})$.

Let $I$ be a maximal independent set of vertices of~$C_5$.  Removing
it from $C_5$ leaves us with an edge and an isolated vertex.  The
colourings of the edge give rise to an $n$-sphere, since $\Hom(K_2,
K_{n+2})\homeo\Sphere^n$.  Similarly the colourings of the isolated
vertex together with the information that this vertex is not isolated
in $C_5$ give rise to another $n$-sphere.  It follows that
$\Hom_I(C_5, K_{n+2})$ can be identified with a subspace of
$\Sphere^n\times\Sphere^n$.  It is not the full space
$\Sphere^n\times\Sphere^n$, since not every multi-colouring of the
edge and the isolated vertex can be extended to~$C_5$; the edge and
the isolated vertex cannot be regarded seperately but interact via the
vertices in~$I$.  It turns out that $\Hom_I(C_5, K_{n+2})$ is the
common boundary of regular neighbourhoods of two subspaces
of~$\Sphere^n\times\Sphere^n$, each corresponding to an element
of~$I$.  These subspaces are ambiently isotopic to the subspaces
$\set{(x,x)\colon x\in\Sphere^n}$ and $\set{(x,-x)\colon
  x\in\Sphere^n}$, which determines $\Hom_I(C_5, K_{n+2})$ up to
homeomorphism and proves the Theorem.  The relevant subdivisions of
the spheres and constructions of regular neighbourhoods are examples
of general constructions that can be carried out for any
pl-triangulation of a manifold, the triangulation being the boundary
of an $(n+1)$-simplex in the case of $\Hom_I(C_5, K_{n+2})$.  These
constructions are carried out in \prettyref{sec:Int}.

\prettyref{sec:C5} puts together the results from the two preceding
sections.  In it we also consider Csorba's stronger conjecture, also
stated in \cite{csorba-thesis}, that the homeomorphism
of~\prettyref{eq:main} can be chosen in such a way that it transports
the involution on $\Hom(C_5,K_{n+2})$ that appears in the Theorem of
Babson and Kozlov mentioned above to the involution
$(x,y)\mapsto(x,-y)$ on $\Sphere^n\times\Sphere^n$.  Our current
approach yields only a proof of a homotopy version of this conjecture,
the complete conjecture remains open.

\begin{acknowledgements}
I would like to thank Péter Csorba for comments on a draft of this
work.  These helped to clarify some points and lead to the inclusion
of the part on the involution on $\Hom(C_5,K_{n+2})$.
\end{acknowledgements}

\section{Preliminaries}
\label{sec:prelims}
We collect some definitions, notations and results that we will need.
\subsection{Posets and cellular complexes}
For a partially ordered set, or \emph{poset}, $P$, we denote by
$\Delta P$ its \emph{order complex}, the simplicial complex with
vertex set $P$ that consists of all chains in~$P$.  A monotone (or
antitone) map $f\colon P\to Q$ between posets induces a simplicial map
$\Delta f\colon\Delta P\to\Delta Q$. All cell complexes are assumed to
be regular.  For a cell complex~$C$ we denote its \emph{face poset} by
$FC$ and its underlying space by $\real C$.  Thus $\Delta(FC)$ is the
barycentric subdivision of the complex~$C$ and
$\real{\Delta(FC)}\homeo\real C$.  If $C$ and~$D$ are cell complexes,
then $C\times D$ denotes the cell complex whose face poset is
isomorphic to $FC\times FD$, even if $C$ and~$D$ are simplicial
complexes.  Since most of our simplicial complexes are order
complexes, we have no need for a notation for the simplicial complex
that is their product, $\Delta(P\times Q)$ is the usual simplicial subdivision
of the cell complex $\Delta P\times\Delta Q$.  Our preferred way to
construct homotopies between maps between order complexes is the
following special case of \cite[Prop.~2.1]{segal-classifying-spaces},
which we will use without citing it.

\begin{lem}
Let $P$, $Q$ be posets.  If $f,g\colon P\to Q$ are order
preserving functions such that $f(p)\le g(p)$ for all $p\in P$, then
the maps $\Delta f,\Delta g\colon\Delta P\to\Delta Q$ are homotopic.
\end{lem}

\begin{proof}
The map $H\colon\set{0,1}\times P\to Q$ defined by $H(0,x)\deq
f(x)$, $H(1,x)\deq g(x)$ is order preserving and hence yields the desired
homotopy
\begin{equation*}
I\times \Delta P\homeo\Delta (\set{0,1}\times P)
\xto{\,\Delta H\,}\Delta Q,
\end{equation*}
 where
we view $\set{0,1}$ as a poset.
\end{proof}
For a poset $P$ its \emph{dual poset $\opp P$} is the poset with the
same elements but the order reversed.  When considering $p\in P$ as an
element of $\opp P$ we write it as $\opp p$.  Thus $p\mapsto\opp p$ is
an antitone bijection $P\to\opp P$.

For a cell $c$, its dimension is denoted by $\adeg c$, the cardinality
of a set~$M$ is written as~$\#M$.  When identifying a simplex~$\sigma$
with the set of its vertices we thus have $\adeg\sigma=\#\sigma-1$.

\subsection{Graphs and graph complexes}
All Graphs that we consider are finite, simple, and without loops.
The vertex set of a graph~$G$ is denoted by $V(G)$, and for $S\subset
V(G)$ the set of all common neighbours of the elements of~$S$ is
denoted by~$\nu(S)$.  A graph homomorphism from $G$ to $H$ is a
function $f\colon V(G)\to V(H)$ that respects the edge relation, i.e.\ such
that $f[\nu(\set v)]\subset\nu(\set{f(v)})$ for all $v\in V(G)$.

\begin{defn}\label{def:ind}
For a graph $G$ we denote by $\ind(G)$ the poset of all independent
subsets of $V(G)$, including the empty set, ordered by inclusion.  By
$\Ind(G)$ we denote the \emph{independence complex of $G$}, i.e.\ the
simplicial complex with vertex set~$V(G)$ and simplices the
independent subsets of~$V(G)$.
\end{defn}

\begin{defn}
Let $G$, $H$ be graphs.  A \emph{multi-homomorphism} from $G$ to~$H$
is a function $\phi\colon V(G)\to\power(V(H))\wo\set\emptyset$ such
that every function $f\colon V(G)\to V(H)$ with $f(v)\in\phi(v)$ for
all~$v\in V(G)$ is a graph homomorphism.
\end{defn}

\begin{defn}\label{def:Hom}
Let $G$, $H$ be graphs.  A function $\phi\colon
V(G)\to\power(V(H))\wo\set\emptyset$ can be identified with a cell of
the cell complex
\[\prod_{v\in V(G)}\Delta^{\#V(H)-1}.\]
The subcomplex of all cells indexed by multi-homomorphisms is denoted
by $\Hom(G, H)$.  We identify elements of $F\Hom(G, H)$ with the
corresponding multi-homomorphisms.  If $f\colon G'\to G$, $g\colon
H\to H'$ are graph homomorphisms, then there is an induced monotone
map
\begin{align*}
F\Hom(f, g)\colon F\Hom(G, H)&\to F\Hom(G', H'),\\
(F\Hom(f,g))(\phi)(v)&\deq g[\phi(f(v))]
\end{align*}
and hence a continous map
\begin{equation*}
\real{\Hom(f, g)}\to\real{\Hom(G, H)}\to\real{\Hom(G', H')}.
\end{equation*}
\end{defn}
The above constructions are functorial.  Details can be found in
\cite{babson-kozlov-i}.

\begin{nota}
Let $n\in\N$.  $P_n$ denotes that path with $n$~edges on the vertex
set $\set{0,\dots,n}$.  $K_n$ denotes the complete graph on
$n$~vertices with vertex set $\set{1,\dots,n}$.  $C_n$ is the cycle of
length~$n$ with vertex set $\set{1,\dots,n}$.
\end{nota}

\subsection{Piecewise linear topology}

\begin{defn}
Let $Y$ be a simplicial complex and $X$ a subcomplex of $Y$.  The subcomplex
\begin{equation*}
N\deq\set{\sigma\in Y\colon\text{There is a $\tau\in X$ with
    $\sigma\unite\tau\in Y$}}
\end{equation*}
of $Y$ is called the \emph{simplicial neighbourhood of $X$ in~$Y$} and
\begin{equation*}
\dot N\deq\set{\sigma\in N\colon\text{$\sigma\intersect\tau=\emptyset$ for
    all $\tau\in X$}}
\end{equation*}
is called the \emph{simplicial neighbourhood boundary}.
We also call the pair $(N, \dot N)$ the simplicial neighbourhood.
\end{defn}

We will need the following form of the Simplicial Neighbourhood
Theorem \cite[3.11]{rourke-sanderson}.
\begin{thm}\label{thm:snt}
Let $M$ be a pl-triangulated compact manifold, $X$ a full subcomplex
of~$M$, and $(N, \dot N)$ the simplicial neighbourhood of $X$ in~$M$.
If $(N, \dot N)$ is a manifold with boundary, then $N$ is a regular
neighbourhood of $X$ in~$M$.
\end{thm}

\section{Restriction maps of \texorpdfstring{$\Hom$}{Hom}-complexes}
\label{sec:res}
Given graphs $G$ and $H$ and an independent subset $I$ of $V(G)$ there
is a subcomplex of $\Hom(G\wo I, H)$ which is homotopy equivalent to
$\Hom(G,H)$.  The use of this smaller complex in the study of the
homotopy type of $\Hom(G,H)$ has been introduced in
\cite[Sec.~2.8]{csorba-thesis}.  The same idea has been implicitly
used in \cite{codd} to study the complexes $\Hom(C_{2r+1}, K_{n+2})$.

In \prettyref{lem:res-homeo} we give a criterion that allows us to
detect cases in which the two spaces are actually homeomorphic.  We
will use this to study $\Hom(C_5, K_{n+2})$.

We round off the investigation by showing that the applicability of the 
Lemma to~$C_5$ is not incidental but a consequence of the fact that 
$\Hom(C_5, K_n)$ is a manifold.

For the results on $\Hom(C_5, K_{n+2})$, only 
\prettyref{def:res}, 
\prettyref{lem:res-homeo},
and \prettyref{ex:c5} will be needed from this section.

\subsection{Homotopy equivalences}
We recall the definition of $\Hom_S(G, H)$ from \cite{csorba-thesis}.
The proofs and examples given here are for expository purposes only.
Neither are they needed for the homeomorphism results we present
afterwards, nor do we claim any originality here, except perhaps for
the point of view described in \prettyref{ex:P2}.
\begin{propdef}[\cite{csorba-thesis}]\label{prop:res-basics}\label{def:res}
Let $G$, $H$ be graphs and $S\subsetneq V(G)$.  Let $i\colon G\wo S\to G$
denote the inclusion of the subgraph induced on the complement of $S$.
We consider the continuous map
\begin{equation*}
\Hom(i,H)\colon\Hom(G,H)\to\Hom(G\wo S,H)
\end{equation*}
and denote its image by $\Hom_S(G, H)$.
\begin{enumerate}
\item\label{it:b1} The map $\Hom(i,H)$ is cellular and 
  $\Hom_S(G, H)$ is a subcomplex of the cell-complex $\Hom(G\wo S, H)$.
\item\label{it:b2} If $S$ is independent, then 
  $\Hom_S(G,H)$ consists of all cells of $\Hom(G\wo S, H)$ 
  indexed by multihomomorphisms
  $\phi$ with $\nu\left(\Union_{u\in\nu(\set v)}\phi(u)\right)\ne\emptyset$
  for all $v\in S$.
\item\label{it:b3} If $S$ is independent, then the map
  \[\Hom(G, H)\to\Hom_S(G, H)\]
  induced by $\Hom(i,H)$ is a homotopy equivalence.
\end{enumerate}
\end{propdef}

\begin{proof}
The first two claims are obvious.  
Let us denote the inclusion $S\to G$ by~$j$.  Then
\begin{equation*}
(F\Hom(i,H),F\Hom(j,H))\colon F\Hom(G,H)\to F\Hom(G\wo S, H)\times F\Hom(S, H)
\end{equation*}
is an embedding.  If $S$ is independent, then the image of this
embedding is \[\set{(\phi,\rho)\in\im F\Hom(i,H)\times F\Hom(S,
  H)\colon \rho\le a(\phi)}\] where $a$ is the antitone map
$\im F\Hom(i,H)\to F\Hom(S, H)$ given by \[a(\phi)(v)\deq
\nu\left(\Union_{u\in\nu(\set v)}\phi(u)\right).\]  The result follows
from \prettyref{lem:aQ}
\end{proof}

We illustrate the Proposition by deriving two well-known facts from
it.

\begin{example}\label{ex:tree-H}
Let $H$ be a graph, $T$ be a tree with at least three vertices and
$v$~a leaf of~$T$.  Since any multihomomorphism from $T\wo\set v$ can
be extended to one from $T$ by mapping $v$ to one of the vertices
assigned to one of the neighbours in $T\wo\set v$ of the neighbour of
$v$, the restriction map $\Hom(T,H)\to \Hom(T\wo\set v,H)$ is
surjective and hence a homotopy equivalence.  By induction any
inclusion $i\colon K_2\to T$ of a an edge in a tree induces a homotopy
equivalence $\Hom(T, H)\homot\Hom(K_2, H)$.  This example generalizes
to folds in the first parameter of~$\Hom$, see~\cite{kozlov-folds} for
this concept.
\end{example}

\begin{example}\label{ex:K2-H}
$\Hom(K_2, H)$ is one of the constructions called the \emph{box
    complex} of $H$, while $\Hom(K_1, H)$ is just the full simplex on
  the vertices of~$H$.  Considering an inclusion $i\colon K_1\to K_2$,
  the image of the map $\Hom(i,H)$ consists of the simplices
  corresponding to sets $A\subset V(H)$ with $\nu(A)\ne\emptyset$.
  This is the \emph{neighbourhood complex} of~$H$ introduced in
  \cite{lovasz}.  Thus the box complex is homotopy equivalent to the
  neighbourhood complex.  This example also appears
  in~\cite{csorba-thesis}.  For comparison of several graph complexes
  also see \cite{matousek-ziegler,shore-subdiv,csorba-box,wiposets}.
\end{example}

We look more closely at $\Hom(P_2, K_{n+2})$.  This will give us an
opportunity to introduce some ideas that can be generalised, but are
not needed for the remainder of this article.

\begin{example}\label{ex:P2}
The relevant inclusions related to maximal independent subsets of $P_2$ are
$j\colon K_1\isom\set1\to P_2$ and $k\colon
P_2\wo\set1\to P_2$.  For $M\subset V(K_{n+2})$, $M\ne\emptyset$, we
define $X_M\deq\set{\phi\in F\Hom(P_2, K_{n+2})\colon \phi(1)=M}$ and
$Y_M\deq\set{\phi\cmps k\colon \phi\in X_M}$.  The sets
$\Delta(Y_{\set u})$, $u\in V(K_{n+2})$ cover $\im\Hom(k, K_{n+2})$,
and $Y_M\intersect Y_N=Y_{M\unite N}$ for all~$M,N$.  Also
$Y_{V(K_{n+2})}=\emptyset$. So, if we define the poset
$Q\deq\power(V(K_{n+2}))\wo\set{\emptyset,V(K_{n+2})}$ ordered by
reverse inclusion, then $Y$ is a diagram of spaces over $Q$, where the
maps are given by inclusion, and $\Hom_{\set 1}(P_2,K_{n+2})\homeo\colim Y$.  
\begin{figure}
\begin{equation*}
\xymatrix{
\set1&\set2&\set3\\
\set{1,2}\ar[u]\ar[ru]&\set{1,3}\ar[lu]\ar[ru]&\set{2,3}\ar[lu]\ar[u]
}
\end{equation*}
\caption{\label{fig:posetQ}}
\end{figure}
\begin{figure}
\input{diag.pstex_t}
\caption{\label{fig:diag}}
\end{figure}
\begin{figure}
\input{colim.pstex_t}
\caption{\label{fig:colim}$\Hom_{\set 1}(P_2,K_3)$}
\end{figure}
\begin{figure}
\input{hcolim.pstex_t}
\caption{\label{fig:hcolim}$\Hom(P_2,K_3)$}
\end{figure}
For the example~$n=1$ the poset~$Q$ and the diagram~$Y$ are shown in
\prettyref{fig:posetQ} and \prettyref{fig:diag}, and $\Hom_{\set
  1}(P_2,K_3)$ is depicted in \prettyref{fig:colim}.  On the other
hand $\Hom(P_2,K_{n+2})$, which is shown in \prettyref{fig:hcolim},
can be seen to be homeomorphic to $\hcolim Y$, and $\Hom(k,K_{n+2})$
corresponds to the natural map from $\hcolim Y$ to $\colim Y$ which is
a homotopy equivalence by the Projection Lemma, since $Y$ is a free diagram  
because of its construction from a cover.  This
indicates another way in which \prettyref{prop:res-basics} could be
proved. We also see that in general $\Hom(G,H)$ will not be
homeomorphic to $\Hom_I(G, H)$, even if $I$~is independent.

The diagram $Y$ can also be used to describe the map $\Hom(j,K_{n+2})$
as the map from the homotopy colimit of~$Y$ to the homotopy colimit of
the constant diagram over $Q$ which assigns a point to every element.
The homotopy colimit of this constant diagram is just~$\Delta Q$,
which is homeomorphic to~$\Sphere^n$.  This map is a homotopy
equivalence, because all $Y_M$ are contractible, which can be seen to
be a consequence of the independence of the set $\set{0,2}$.  This
description is nearer to the proof of \prettyref{prop:res-basics}
given above.  In general the poset $Q$ would have to be replaced by
the face poset of $\Hom_S(G,H)$.  The situation in this example is
special, because $P_2$ is bipartite.
\end{example}

\subsection{Homeomorphisms}
We examine situations in which the complexes $\Hom(G, K_n)$ and
$\Hom_I(G, K_n)$ are indeed homeomorphic.  The independence complex
of~$G$ (\prettyref{def:ind}) plays an important role in recognising
these situations.
\begin{lem}\label{lem:res-homeo}
Let $G$ be a graph, $n\in\N$, and $S\subsetneq V(G)$ an independent set.
We define subsets $A_v$
of~$\ind(G)$ for all~$v\in V(G)$ by
\begin{align*}
A_v&\deq\set{I\in\ind(G)\colon v\in I}
\intertext{and similarly}
B_v&\deq\begin{cases}
\set{I\in\ind(G\wo S)\colon v\in I},&v\notin S\\
\set{I\in\ind(G\wo S)\colon I\unite\set v\in\ind(G)},&v\in S.
\end{cases}
\end{align*}
If there is a homeomorphism $h\colon
\abs{\Delta(\ind(G))}\to\abs{\Delta(\ind(G\wo S))}$ with
$h[\abs{\Delta(A_v)}]=\abs{\Delta(B_v)}$ for all~$v\in V(G)$, then
$\Hom(G,K_n)\homeo\Hom_S(G,K_n)$.
\end{lem}

\begin{proof}
We look at $\Hom(G,K_n)$ one colour at a time instead of one vertex at
a time.  The map
\begin{align*}
\gamma\colon F\Hom(G,K_n)&\to \prod_{c\in V(K_n)}\ind(G)\\
\phi&\mapsto(\set{u\in V(G)\colon c\in\phi(u)})_c
\end{align*}
is a well-defined poset embedding.  Hence $\Hom(G,K_n)$ is
homeomorphic to the order complex of the image of~$\gamma$.  The poset
$\ind(G)^{V(K_n)}$ describes multicolourings without the condition
that each vertex has to obtain at least one colour.  Spelling out this
additional condition yields
\begin{equation*}
\image\gamma=\Intersection_{v\in V(G)}\Union_{d\in V(K_n)}
 \prod_{c\in V(K_n)}\begin{cases}
                       A_v,& c=d\\
                       \ind(G),& c\ne d.
 \end{cases}
\end{equation*}
Since the $A_v$ are closed from above, the operation of taking the
order complex commutes with all constructions appearing in this
equation, and
\begin{equation*}
\Hom(G, K_n)\homeo\Intersection_{v\in V(G)}\Union_{d\in V(K_n)}
 \prod_{c\in V(K_n)}\begin{cases}
                       \Delta(A_v),& c=d\\
                       \Delta(\ind(G)),& c\ne d.
 \end{cases}
\end{equation*}
follows.  In the same way $F\Hom(G\wo S, K_n)$ embeds in $\ind(G\wo
S)^{V(K_n)}$.  $\Hom_S(G, K_n)$ is distinguished in this
embedding by the condition that for every vertex in $S$ one of the
colours is not used by any of its neighbours.  Since $B_v$ is closed
from above for $v\in S$ and closed from below for $v\notin S$, the
above argument can be repeated to obtain
\begin{equation*}
\Hom_S(G, K_n)\homeo\Intersection_{v\in V(G)}\Union_{d\in V(K_n)}
 \prod_{c\in V(K_n)}\begin{cases}
                       \Delta(B_v),& c=d\\
                       \Delta(\ind(G\wo S)),& c\ne d.
 \end{cases}
\end{equation*}
It is now obvious how the homeomorphism~$h$ induces a homeomorphism
$\Hom(G, K_n)\homeo\Hom_S(G, K_n)$.
\end{proof}

Continuing \prettyref{ex:K2-H} we give a quite trivial application of
the Lemma.
\begin{example}
We consider the inclusion $i\colon K_1\to K_2$.
$\ind(K_2)=\set{\emptyset,\set0,\set1}$, hence $\Delta(\ind(K_2))$ is
homeomorphic to an interval.  The two points of the boundary are
$\Delta(\set{\set0})=\Delta(A_0)$ and
$\Delta(\set{\set1})=\Delta(A_1)$.  On the other hand
$\ind(K_1)=\set{\emptyset,\set0}$, so $\Delta(\ind(K_1))$ is also
homeomorphic to an interval.  The boundary points are
$\Delta(\set{\set0})=\Delta(B_0)$ and
$\Delta(\set{\set1})=\Delta(B_1)$.  Thus the Lemma is applicable, and
$\Hom(K_2, K_{n+2})\homeo\Hom_{\set 1}(K_2, K_{n+2})$.  The latter complex
is the boundary of an $(n+1)$-simplex.  This gives yet another proof
of $\Hom(K_2, K_{n+2})\homeo\Sphere^n$.
\end{example}

\begin{example}
The Lemma cannot be applied to the restriction maps considered in
\prettyref{ex:P2} as we have seen there.  Indeed $\Delta(\ind(P_2))$
is not homeomorphic to $\Delta(\ind(K_1))$, $\Delta(\ind(K_2))$, or
$\Delta(\ind(P_2\wo\set1))$.  $\ind(P_2)$ has the maximal elements
$\set1$ and $\set{0,2}$ and hence is not a pure complex.  As we have
seen $\Delta(\ind(K_1))$ and $\Delta(\ind(K_2))$ are intervals.
$\ind(P_2\wo\set1)\isom\ind(K_1)\times\ind(K_1)$ and hence
$\Delta(\ind(P_2\wo\set1))$ is homeomorphic to a $2$-disk.
\end{example}

A less trivial example and our main reason for this investigation.
\begin{example}\label{ex:c5}
We consider $\set{2,4}\subset V(C_5)$.
\prettyref{fig:indc5} shows that \prettyref{lem:res-homeo}
can be applied and hence
$\Hom(C_5, K_{n+2})\homeo\Hom_{\set{2,4}}(C_5, K_{n+2})$.
\begin{figure}
\input{indc5.pstex_t}
\caption{\label{fig:indc5}$\Delta(\ind C_5)$ vs.\ $\Delta(\ind(C_5\wo\set{2,4}))$}
\end{figure}
\end{example}

\subsection{Manifolds}
In the two examples we have just seen where $\Hom(G,
K_n)\homeo\Hom_I(G, K_n)$ for a non-empty independent set~$I$ and
all~$n$ the independence complex~$\Ind(G)$ is a sphere,
$\Ind(K_2)\homeo\Sphere^0$ and $\Ind(C_5)\homeo\Sphere^1$.  This is
what makes $\Hom(G, K_n)$ a manifold in these cases.  Complexes
$\Hom(G, K_n)$ that are manifolds have been investigated
in~\cite{csorba-lutz}.  We show that for these 
$\Hom(G,K_n)\homeo\Hom_I(G, K_n)$ always holds.

\begin{thm}[\cite{csorba-lutz}]\label{thm:hom-mf}
Let $G$ be a graph.  The following statements are equivalent.
\begin{enumerate}
\item$\Ind(G)$ is a pl-sphere.
\item$\Hom(G, K_n)$ is a pl-manifold for all~$n$.
\item$\Hom(G, K_{\chi(G)+1})$ is a pl-manifold.
\end{enumerate}
If these statements hold, then $\dim\Hom(G, K_n)=n(\dim\Ind(G)+1)-\abs{V(G)}$.
\end{thm}
We repeat the proof from \cite{csorba-lutz}, since we will afterwards build upon
the ideas used in it.
\begin{proof}
Let $\phi\in F\Hom(G, K_n)$.  For every $k\in V(K_n)$ the cell~$\phi$
determines an independent set $\set{v\in V(G)\colon k\in\phi(v)}$ and
hence a simplex~$\sigma_k$ of~$\Ind(G)$, where we allow the empty
simplex.  The link of $\phi$ in $\Hom(G,K_n)$ is homeomorphic to the
join of the links of the $\sigma_k$ in $\Ind(G)$.  If
$\Ind(G)\homeo\Sphere^m$, then the link of $\sigma_k$ is an
$(m-1-\abs{\sigma_k})$-sphere, and the link of~$\phi$ is a sphere of
dimension
\[n-1+\sum_k(m-1-\abs{\sigma_k})
 =n(m+1)-1-n-\sum_k\abs{\sigma_k} =n(m+1)-\#(V(G))-1-\abs\phi,\] so
 $\Hom(G,K_n)$ is a manifold of dimension $n(m+1)-\#(V(G))$.  If
 $n>\chi(G)$ then there is a cell~$\phi$ for which the empty simplex
 is among the $\sigma_k$, so the link of $\phi$ is a join of spaces of
 which one is $\Ind(G)$.  For the join to be a sphere it is necessary
 for $\Ind(G)$ to be a sphere, see \cite[2.24(5)]{rourke-sanderson}.
\end{proof}
We will want to apply \prettyref{lem:res-homeo} in this case and
therefore start examining the relationship between $\Ind(G)$
and~$\Ind(G\wo S)$.
\begin{lem}
Let $G$ be a graph and $S\subset V(G)$ a non-empty independent set.
Then $\Ind(G\wo S)$ is the complement of the interior of the
simplicial neighbourhood of the simplex~$S$ in $\Ind(G)$.  If
$\Ind(G)$ is a pl-manifold, then this simplicial neighbourhood is a
regular neighbourhood.  Hence if $\Ind(G)\homeo\Sphere^n$ then
$\Ind(G\wo S)\homeo\Disk^n$.
\end{lem}

\begin{proof}
We only have to show that the simplicial neighbourhood of $S$ is a
regular neighbourhood if $\Ind(G)$ is a manifold.  We assume that
$\Ind(G)$ is an $n$-manifold. The simplicial neighbourhood of the
simplex~$S$ consists of all $I\in\ind(G)\wo\set\emptyset$ such that
there exists a $v\in S$ such that $I\unite\set v\in\Ind(G)$.  Its
boundary consists of those~$I$ for which additionally $I\intersect
S=\emptyset$.  By \prettyref{thm:snt} it will be
sufficient to show that the simplicial neighbourhood is a manifold
with boundary, the boundary being as just described.  To prove this we
take a simplex~$I$ in the boundary and examine its link.  Its link
in~$\Ind(G)$ can be identified with~$\Ind(G')$ where $G'$ is the graph
obtained from~$G$ by deleting the vertices in~$I$ and all of their
neighbours and is an $(n-\#I)$-sphere, since $\Ind(G)$ is an
$n$-manifold.  We set $S'\deq S\intersect V(G')$. Since $I$ is a
simplex of the simplicial neighbourhood of $S$, $S'$~is non-empty.  In
the link of~$I$ the boundary of the simplicial neighbourhood of~$S$
corresponds to the simplicial neighbourhood of~$S'$ in $\Ind(G')$,
which is an unlinked $(n-1-\#I)$-sphere by induction.
\end{proof}

\begin{rem}
Since $\Ind(G)$ determines~$G$, it is easy to describe, which
simplicial complexes are independence complexes of graphs.  These are
the flag simplicial complexes.  The preceding Lemma can therefore be
viewed as a statement on flag complexes.
\end{rem}

\begin{lem}\label{lem:mf-homeo}
Let $G$ be a graph, $n\in\N$, and $S\subsetneq V(G)$ an independent
set.  If $\Ind(G)$ is a pl-sphere, then
$\Hom(G,K_n)\homeo\Hom_S(G,K_n)$.
\end{lem}

\begin{proof}
We want to apply \prettyref{lem:res-homeo} and use the notation
introduced there.  Let $m\deq\dim\Ind(G)$.  $\Delta(\ind(G))$ is a
cone over $\Ind(G)$, therefore an $(m+1)$-ball, and $A_v$ is the star
of~$v$ in the barycentric subdivision of~$\Ind(G)$, i.e. the cell of
the dual complex which is dual to $v$.  We will now describe
$\Delta(\ind(G\wo S))$.

If $S\ne\emptyset$, then we have seen in the preceding Lemma that
$\Ind(G\wo S)$ is an $m$-ball.  This ball is covered by the sets
$\Delta(B_v)$ with $v\in V(G)\wo S$.  The boundary of this ball
consists of the simplices $I$ for which there exists a $v\in S$ such
that $I\unite\set v$ is independent.  Therefore the sets $\Delta(B_v)$
with $v\in S$ cover the cone over this boundary in $\Delta(\ind(G\wo
S))$.  It follows that $\Delta(\ind(G\wo S))$ is an $(m+1)$-ball and
that the sets $\Delta(B_v)$ for all $v\in V(G)$ cover the boundary of
this ball.  The last statement is also true if $S=\emptyset$, and we
will from here on allow this case in order to facilitate inductive
arguments.

For $M\in\ind(G)\wo\set\emptyset$ we set
\begin{equation*}
C_M\deq\set{I\in\ind(G\wo S)\colon
  \text{$M\wo S\subset I$, $I\unite M\in\ind(G)$}}.
\end{equation*}
Then $B_v=C_{\set v}$, and we have $\Delta(C_M)\intersect
\Delta(C_N)=\Delta(C_{M\unite N})$, if $M\unite N$ is independent, and
$C_M\intersect C_N=\emptyset$ otherwise.

Let $M\in\ind(G)\wo\set\emptyset$ and let $G'$ be the graph obtained
from $G$ by removing all vertices in~$M$ and their neighbours.  As
seen in the proof of the preceding Lemma, $\Ind(G')$ is homeomorphic
to the link of $M$ in $\Ind(G)$ and hence an $(m-\#M)$-sphere.  Let
$S'\deq S\intersect V(G')$.  The map
\begin{align*}
C_M&\to\ind(G'\wo S')\\
I&\mapsto I\wo M
\end{align*}
is an isomorphism and for $N\supsetneq M$ it maps $C_N$ to
\[\set{I\in\ind(G'\wo S')\colon
  \text{$(N\wo M)\wo S\subset I$, $I\unite (N\wo M)\in\ind(G')$}}.\]
Therefore by induction $\Delta(C_M)$ is an $(m-\adeg M)$-ball and its
boundary is the union of the $\Delta(C_N)$ for $N\supsetneq M$.

This shows that the sets $\Delta(C_M)$ form a cell-decomposition of
the boundary of $\Delta(\ind(G\wo S))$ that has a face poset that is
isomorphic to the face poset of the dual complex of the simplicial
complex $\Ind(G)$.  This yields a homeomorphism from the boundary of
the ball $\Delta(\ind(G))$ to the boundary of the ball
$\Delta(\ind(G\wo S))$ that carries $\Delta(A_v)$ to $\Delta(B_v)$ for
every $v\in V(G)$.  This homeomorphism can be extended to a homeomorphism
between balls and \prettyref{lem:res-homeo} can be applied.
\end{proof}

\section{Subdivisions and diagonal approximations}
\label{sec:Int}
We examine several constructions that occur naturally in the study of
homomorphism complexes of graphs, although they will not be mentioned
in this section. All posets in this section are assumed to be finite.

\subsection{Edge subdivision}
For an ordered simplicial complex there is a subdivision that
introduces a vertex for every edge, which is why we will call it its
\emph{edge subdivision}.  For an order complex of a poset, it is the
order complex of the interval poset.  It has been first described in
\cite{walker-posets}.  In connection with graph complexes it has been
used in \cite{wiposets}.
\begin{defprop}
For a poset $P$ we define the \emph{closed interval poset of $P$} to
be the subposet of $P\times\opp P$ consisting of all $(p,\opp q)$ with
$p\le q$ and denote it by $\Int P$.  $\Int$ is a functor from posets
to posets.  $\Delta(\Int P)$ is a subdivision of $\Delta P$.  More
precisely, choosing a point $x_{p,q}$ in the interior of the simplex
$\simp{p,q}$ of $\Delta P$ for every pair $(p,q)$ with $p<q$
determines a unique map $\real{\Delta(\Int P)}\to\real{\Delta P}$
which sends $(p,\opp p)$ to $p$, $(p,\opp q)$ to $x_{p,q}$ for $p<q$,
and is affine on simplices, and this map is a homeomorphism.  In
particular for $0<\lambda<1$ we define
$h_P^\lambda\colon\real{\Delta\Int P}\to\real{\Delta P}$ by choosing
$x_{p,q}$ as $\lambda p+(1-\lambda)q$.  $h^\lambda$ is a natural
equivalence from $P\mapsto\real{\Delta\Int P}$ to
$P\mapsto\real{\Delta P}$.  We abbreviate $h^{1/2}_P$ as $h_P$.  All
choices of $x_{p,q}$ yield isotopic maps, and any such map is
homotopic to each of the two maps $\Delta(\Int P)\to\Delta P$ given by
$(p,\opp q)\mapsto p$ or $(p,\opp q)\mapsto q$.
\end{defprop}

\begin{proof}
To see that $\Delta(\Int P)$ is a subdivision of~$\Delta P$ one either
checks directly that the map $\real{\Delta(\Int P)}\to\real{\Delta P}$
is bijective as in \cite{walker-posets}, or one uses that for a
simplex with vertices $p_0<p_1<\dots<p_r$ the poset
$\Int\set{p_0,\dots,p_r}$ has the minimum $(p_0,\opp{p_r})$ and the
elements $(p_0,\opp q)$ and $(q',\opp{p_r})$ are comparable only if
one of them equals this minimum.  Hence
$\Delta(\Int\set{p_0,\dots,p_r})$ is a cone over
$\Delta(\Int\set{p_0,\dots,p_{r-1}})\unite\Delta(\Int\set{p_1,\dots,p_r})$,
which yields a recursive description of the subdivision.  To connect
subdivision maps for different choices of $x_{p,q}$ by an isotopy, one
moves the points $x_{p,q}$.  The homotopies are constructed by
extending the definition of $h^\lambda_P$ to $\lambda=0$ and
$\lambda=1$.
\end{proof}

\begin{rem}
$(p_0,\opp{q_0})\le (p_1,\opp{q_1})\iff [p_0,q_0]\supset[p_1,q_1]$,
  hence the name interval poset.
\end{rem}

For illustration we use this subdivision to prove a simple Lemma,
which is quite specialised but useful in the proof of
\prettyref{prop:res-basics}.

\begin{lem}\label{lem:aQ}
Let $P,Q$ be posets, $a\colon P\to Q$ an antitone map, and
\[R\deq\set{(p,q)\in P\times Q\colon q\le a(p)}.\]  Assume that any two
elements of $Q$ having a common upper bound have a unique least upper
bound.  Then the map $\Delta(R)\to\Delta(P)$ induced by projection
onto the first factor is a homotopy equivalence.
\end{lem}

\begin{proof}
We define monotone maps $\pi\colon R\to P$, $(p,q)\mapsto p$ and
$\sigma\colon \Int P\to R$, $(p_0,\opp{p_1}) \mapsto(p_0,a(p_1))$.
Then $(\pi\cmps\sigma)(p_0,\opp{p_1})=p_0$, so
$\real{\pi}\cmps\real\sigma\homot h_P$.  Also
\[(\sigma\cmps\Int\pi)((p_0,q_0),(\opp p_1,\opp
q_1))=(p_0,a(p_1))\le(p_0,a(p_1)\lmax q_0)\ge(p_0,q_0)\] and, since
$((p_0,q_0),(\opp p_1,\opp q_1))\mapsto(p_0,a(p_1)\lmax q_0)$ is a
well-defined monotone map, $\real\sigma\cmps
h_P^{-1}\cmps\real\pi\cmps h_R\homot h_R$.  Thus $\real\sigma\cmps
h_P^{-1}$ is homotopy inverse to $\real\pi$.
\end{proof}

\begin{rem}
In the proof the chain poset of~$P$ would have worked equally well,
but the interval poset seems more natural.  Also we will see that
iterated interval posets are easier to describe than iterated chain
posets.
\end{rem}

\subsection{Diagonal approximation}
By construction the edge subdivision of $\Delta P$ is a subcomplex of
$\Delta(P\times\opp P)\homeo\Delta P\times\Delta P$.  It can be viewed
as a substitute for the diagonal, that is not a subcomplex of
$\Delta(P\times\opp P)$.
\begin{prop}\label{prop:isotopy}
Let $P$ be a poset and $o\colon\Delta P\to\Delta\opp P$ denote the
homeomorphism induced by $p\mapsto\opp p$.  Then the diagonal map
$\Delta(\Int P)
 \xto{h_P}\Delta P
 \xto{(\id, \Delta o)}\Delta P\times \Delta\opp P$ is isotopic to the map
$\Delta(\Int P)\to\Delta(P\times\opp P)\homeo\Delta P\times\Delta\opp P$
induced by inclusion.
\end{prop}

\begin{proof}
For $0\le\lambda\le1$ we consider the map $(h^\lambda_P, \Delta o\cmps
h^{1-\lambda}_P)\colon \Delta(\Int P)\to\Delta P\times\Delta\opp P$.
For $\lambda=1$ this map agrees with the map induced by inclusion and
is hence an embedding.  For $0<\lambda<1$ both components are
homeomorphism, so it is also an embedding.  For $\lambda=\frac12$ the
map agrees with the diagonal map.
\end{proof}

\subsection{Iterated and generalised edge subdivisions}
Iterating the intervall poset construction yields again easy to
describe posets.
\begin{prop}
Let $P$ be a poset.  The map
\begin{align*}
\Int(\Int P)&
\to\set{(p,\opp q,r,\opp s)\in P\times\opp P\times P\times\opp P
    \colon p\le q\le r\le s}\\
\left((x,\opp y),\opp{(u,\opp v)}\right)&
\mapsto(x,\opp u,v,\opp y)
\end{align*}
is an isomorphism of posets.\qed
\end{prop}
This isomorphism makes it natural to expect that $\Delta\set{(p,\opp
  q,r)\in P\times\opp P\times P\colon p\le q\le r}$ would also be a
subdivision of~$\Delta P$, and we will confirm this in the next
Proposition.  It also indicates another connection between
$\Delta(\Int(\Int P))$ and $\Delta(\Int P)$ than the first being the
edge subdivision of the second.  We therefore take a closer look at
the complex $\Delta\Int P$.

An $r$-simplex of $\Delta\Int P$ is of the form
$(p_0,\opp{q_0})<\cdots<(p_r,\opp{q_r})$.  This implies $p_0\le
p_1\le\cdots\le p_r\le q_r\le q_{r-1}\le\cdots\le q_0$.  Hence the
simplex is contained in $\Delta(\set{p_i}\times\{\opp{q_j}\})$ and
$\max\set{p_i}\le\min\{q_j\}$.  On the other hand for chains
$p_0<\dots<p_r$ and $q_0<\dots<q_s$ with $p_r\le q_0$ the subset
$\{p_i\}\times\{\opp{q_i}\}$ of~$P\times\opp P$ is contained in $\Int
P$.  $\Delta(\{p_i\}\times\{\opp{q_i}\})$ is the product of an
$r$-simplex and an $s$-simplex and its image under
$h_P\colon\real{\Delta\Int P}\to\real{\Delta P}$ is contained in the
simplex $\Delta(\{p_i\}\unite\{q_j\})$.  Thus the sub-complex of the
cell-complex $\Delta P\ctimes\Delta P$ consisting of all cells
$\simp{p_0,\dots,p_r}\times\simp{q_0,\dots,q_s}$ with
$\max\set{p_i}\le\min\{q_j\}$ is isomomorphic to a cell-subdivision of
$\Delta P$, and it has $\Delta\Int P$ as a simplicial subdivision.  It
follows that any subdivision of $\Delta P\ctimes\Delta P$ leads to a
subdivision of $\Delta P$. In particular $\Delta\Int P$ arises in this
way by subdividing $\Delta P\ctimes\Delta P$ as $\Delta(P\times\opp
P)$, and $\Delta(\Int(\Int P))$ arises by subdividing $\Delta
P\times\Delta P$ as $\Delta(\Int P\times\Int P)$.  On the other hand,
the subdivision of $\Delta P$ obtained in this way by triangulating
$\Delta P\times\Delta P$ as $\Delta(P\times P)$ is not a full
subcomplex of $\Delta(P\times P)$, which underlines the usefulness of
the interval poset.

\begin{prop}\label{prop:int32}
Let $P$ be a poset.  There is a commutative diagram
\begin{equation*}
\xymatrix{
\real{\Delta\Int P}
\ar[r]
&
\real{\Delta(P\times\opp P)}
\\
\real{\Delta\set{(p,\opp q,r)\in P\times\opp P\times P\colon p\le q\le r}}
\ar[r]
\ar[u]_\homeo
&
\real{\Delta(\Int P\times P)}
\ar[u]_\homeo
}
\end{equation*}
with the horizontal arrows induced by inclusion and the vertical
arrows homeomorphisms.
\end{prop}

\begin{proof}
$\Delta(\Int P\times P)$ is a subdivision of $\Delta P\ctimes\Delta
  P$.  According to the preceding discussion we only have to check
  that $\Delta\set{(p,\opp q,r)\in\Int P\times P\colon q\le r}$ is the
  correct subcomplex.  Let
  $(p_0,\opp{q_0},r_0)<\dots<(p_k,\opp{q_k},r_k)$ be a chain in $\Int
  P\times P$.  The smallest cell of $\Delta P\ctimes\Delta P$ that
  contains it is $\Delta(\{p_i\}\unite\{q_i\})\times\Delta\set{r_i}$.
  Now $\max(\{p_i\}\unite\{q_i\})=q_0$, $\min\{r_i\}=r_0$, and $q_0\le
  r_0$ if and only if $q_i=r_i$ for all~$i$.
\end{proof}

\begin{rem}
While $\Delta\Int P$ is a subdivision of $\Delta P$ and hence,
as we have used, $\Delta(\Int P\times P)$ and $\Delta(P\times\opp
P)$ are both subdivisions of $\Delta P\ctimes\Delta P$,
the complex $\Delta(\Int P\times P)$ is not a subdivision of
$\Delta(P\times\opp P)$.
\end{rem}

\subsection{Manifolds}
In general there is no ambient isotopy from the inclusion map
$\Delta(\Int P)\to\Delta(P\times\opp P)$ to the diagonal map, as the
case where $\Delta P$ is a simplex shows.  If $\Delta P$ is a manifold
however, then there is one.
\begin{prop}\label{prop:ambient-isotopy}
Let $P$ be a poset such that $\Delta P$ is a pl-triangulation of a compact 
manifold $M^n$.
Then the isotopy from \prettyref{prop:isotopy} is ambient.
\end{prop}

The technical part of the following proof is needed only for $n=2$.
In \prettyref{thm:result} this case is used for $\Hom(C_5,
K_4)\homeo\R P^3$, which has already been proved in
\cite{csorba-lutz}.

\begin{proof}
The codimension of $M$ in $M\times M$ equals~$n$.  Since any isotopy
in codimension at least~$3$ is ambient (\cite{unknot}, see also
\cite[7.3]{rourke-sanderson}) and a simple drawing (\prettyref{fig:ambient}) 
should convince the
reader of the truth of the Proposition for~$n=1$, we will spell out
the details only for~$n=2$.
\begin{figure}
\input{ambient.pstex_t}
\caption{\label{fig:ambient}}
\end{figure}

Let $n=2$.  We will construct a sequence of subpolyhedra of $M\times
M$, the first one being the diagonal, the last one being $\Delta(\Int
P)$, such that for any consecutive pair there is an isotopy of
$M\times M$ which is equal to the identity outside a small
neighbourhood of a $4$-cell of the form $\sigma\times\sigma$ with
$\sigma$ a simplex of $\Delta P$ and carries one polyhedron into the
other.  Each of these subpolyhedra will be of the following form.  For
$p_0<p_1<p_2$, $p_i\in P$, it will contain either the diagonal of
$\simp{p_0,p_1,p_2}\times\simp{p_0,p_1,p_2}$ or
$\Delta(\Int\set{p_0,p_1,p_2})$, i.e.  the union of the simplices
$\simp{p_0}\times\simp{p_0,p_1,p_2}$,
$\simp{p_0,p_1}\times\simp{p_1,p_2}$, and
$\simp{p_0,p_1,p_2}\times\simp{p_2}$.  Additionally, for $p_0<p_1$
such that for exactly one of the two $2$-simplices adjacent to
$\simp{p_0,p_1}$ the polyhedron contains the corresponding part of the
diagonal, the simplex $\simp{(p_0,p_0),(p_0,p_1),(p_1,p_1)}$ will be
included.  Now assume we are given two such subpolyhedra $Q_0$, $Q_1$
which differ only for the choice for the simplex $\simp{p_0,p_1,p_2}$.
We assume that $Q_0$ agrees on
$D\deq\simp{p_0,p_1,p_2}\times\simp{p_0,p_1,p_2}$ with the diagonal,
and that we know by induction that $Q_0$ is a locally flat submanifold
of $M\times M$.  Let $N$ be a small closed regular neighbourhood of
$D$.  It suffices to show that both $(N,N\intersect Q_i)$ are
unknotted ball pairs, because there is then an isotopy of $N$ which is
constant on the boundary of $N$ and moves $N\times Q_0$ to $N\times
Q_1$ \cite[4.4, 3.22(i)]{rourke-sanderson}, and that can be extended to
$M\times M$ by the identity.  Since $Q_0$ is a submanifold and we will
see that $D\intersect Q_0$ is a $2$-ball with boundary
$\closure{Q_0\wo D}$, we know that if $N$ is chosen small enough, then
$(N,N\intersect Q_i)$ is homeomorphic to $(D,D\intersect Q_i)$ with an
external collar attached to $(\dd D, \closure{Q_0\wo D})$.  Since
$Q_i$ is contained in the union of all cells of the form
$\sigma\times\sigma$ with $\sigma$ a simplex of $\Delta P$ and the
intersection of two of these is either empty or also of this form,
$\closure{Q_0\wo D}$ is the union of either
$\simp{p_i}\times\simp{p_i,p_j}\unite\simp{p_i,p_j}\times\simp{p_j}$
or $\simp{(p_i,p_i),(p_j,p_j)} $for every pair $i<j$ .  
$D\intersect Q_i$ is a two ball with
boundary $\closure{Q_0\wo D}$.  $D\intersect Q_0$ consists of
$\simp{(p_0,p_0),(p_1,p_1),(p_2,p_2)}$ plus simplices
$\simp{(p_i,p_i),(p_i,p_j),(p_j,p_j)}$ as needed.  $D\intersect Q_1$
consists of $\Delta(\Int\set{p_0,p_1,p_2})$ plus the remaining
simplices of the form $\simp{(p_i,p_i),(p_i,p_j),(p_j,p_j)}$.  It can
therefore be checked without further knowledge about $M$ that
$(N,N\intersect Q_i)$ is an unknotted $(4,2)$-ball pair for
$i\in\set{0,1}$.
\end{proof}

\begin{figure}
\input{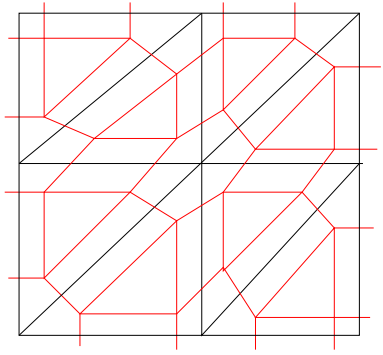}
\caption{\label{fig:IntM}The cell-complex with face poset $\Int M$.}
\end{figure}

If $P$ is the face poset of a manifold~$M$, then there is another
useful way to see $\Delta(\Int P)$, which is the edge subdivision of
the barycentric subdividion of~$M$, as a subdivision of a
cell-complex.  In this case $\opp P$ is the face poset of the dual
cell-complex $\opp M$, and hence $P\times\opp P$ is the face poset of
$M\times\opp M$.  $\Delta(\Int P)$ is then the barycentric subdivision
of a subcomplex of $M\times\opp M$, namely the subcomplex of all cells
$\sigma\times\opp\tau$, where $\tau$ is a simplex of $M$, $\opp\tau$
its dual cell, and $\sigma$ a face of $\tau$.  The maximal cells are
those with $\sigma=\tau$.  Thus $\Delta(\Int P)$ is the barycentric
subdivision of a cell-complex which is homeomorphic to $M$ and
contains a facet for every face of $M$.  
\draftcomment{Add  references?}
\prettyref{fig:IntM} shows this cell-complex for a
part of a manifold.  \prettyref{fig:ambient} shows it as a subcomplex
of $M\times\opp M$, where $M$ is the boundary of a $2$-simplex.

\begin{prop}\label{prop:nb-mf}
Let $P$ be the face poset of a pl-triangulation of a compact manifold $M^n$.
Let $N,B\subset\opp P\times\Int P$ be given by
\begin{align*}
N&\deq\set{(\opp p,q,\opp r)\colon q\le r,p\intersect r\ne\emptyset},\\
B&\deq\set{(\opp p,q,\opp r)\colon q\le r,p\intersect r\ne\emptyset,p\not\le q},
\end{align*}
where notationally identified elements of~$P$ with sets of vertices.
Then $\Delta N$ is a $2n$-dimensional manifold with boundary~$\Delta B$.
\end{prop}

\begin{figure}
\input{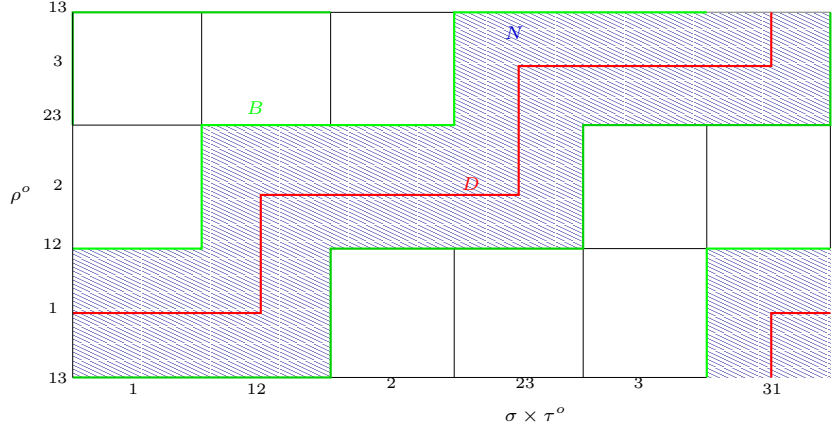}
\caption{The situation of Propositions~\ref{prop:nb-mf} 
and~\ref{prop:nbd} for~$M=\dd\Delta^2$.}
\end{figure}

\begin{proof}
\def\sdiff{\triangle}
\def\void{\emptyset}
\def\emptycmplx{\Sphere^{-1}}
Let $s=(\opp p, q, \opp r)\in\opp P\times\Int P$.  Then 
$\Delta\set{t\in\opp
  P\times\Int P\colon t<s}\homeo \Sphere^{2n-\abs p-\abs r+q-1}$.  So
$\opp P\times\Int P$ is the face poset of a cell complex with $N$
and~$B$ corresponding to subcomplexes, and it will be sufficient to
describe the links of cells of this complex.  We define
\begin{align*}
\lk(s)&\deq\set{t\in\opp P\times\Int P\colon t>s},\\
\lk^+(s)&\deq
    \set{(\opp u, v,\opp w)\in\lk(s)\colon u\intersect w\ne\emptyset},\\
\lk^-(s)&\deq
    \set{(\opp u, v,\opp w)\in\lk(s)\colon u\not\subset v}.
\end{align*}
For $s\in N$ we will have to show that
\begin{equation}\label{eq:a}
(\Delta\lk(s);\Delta\lk^+(s),\lk^-(s))\homeo
\begin{cases}
  (\Sphere^{d-1};\Sphere^{d-1},\emptyset),&s\in N\wo B\\
  (\Sphere^{d-1};\Disk^{d-1}_+,\Disk^{d-1}_-),&s\in B
\end{cases}
\end{equation}
with $d=\abs r-\abs q+\abs p$.  

To be able to approach this inductively, we consider the link of~$s$
in a direction given by a set~$A$ of vertices. 
For $A\subset p\unite r$, we make the following
auxiliary definitions.
\begin{align*}
  \lk_A(s)&\deq\set{(\opp u,v,\opp w)\in\opp{\hat P}\times\Int P\colon
     \text{$(\opp u,v,\opp w)>s$,
     $u\sdiff p\subset A$, 
     $v\sdiff q\subset A$,
     $w\sdiff r\subset A$}},\\
  \lk_A^+(s)&\deq
    \set{(\opp u, v,\opp w)\in\lk_A(s)\colon 
           u\intersect w\intersect A\ne\emptyset},\\
  \lk_A^-(s)&\deq
    \set{(\opp u, v,\opp w)\in\lk_A(s)\colon 
           u\intersect A\not\subset v}.
\end{align*}
Here $\hat P$ is $P$ with a minimum, the empty simplex, added, and
$u\Delta p$ denotes the symmetric difference $(u\unite
p)\wo(u\intersect p)$.

For $A\intersect B=\emptyset$ we have
\begin{equation}\label{eq:b}
\begin{split}
(\Delta\lk_{A\unite B}(s);\Delta\lk_{A\unite
    B}^+(s),\Delta\lk_{A\unite B}^-(s))
  \homeo
&(\Delta\lk_{A}(s)\join\Delta\lk_B(s);\\
&\quad
  \Delta\lk_A^+(s)\join\Delta\lk_B(s)\unite\Delta\lk_A(s)\join\Delta\lk_B^+(s),\\
&\quad
  \Delta\lk_A^-(s)\join\Delta\lk_B(s)\unite\Delta\lk_A(s)\join\Delta\lk_B^-(s)).
\end{split}
\end{equation}
However, for this to be correct, care has to be taken with regard to
empty sets.  We distinguish between the void complex~$\emptyset$ with
$X\join\emptyset=\emptyset$ and the empty complex~$\Sphere^{-1}$ with
$X\join\Sphere^{-1}=X$, which is consistent with $X\join\Sphere^{k-1}$
being the $k$-fold suspension of~$X$.  We set
$\Delta\lk^+_A(s)=\void$, if and only if $A\intersect p\intersect
r=\emptyset$.  The justification for this is that we regard $s$ as the
tip of a cone over $\lk_A(s)$ and want $s$ to be in the cone over
$\lk_A^+(s)$ only if $A\intersect p\intersect r\ne\emptyset$.  In the
same way we set $\lk_A^-(s)=\void$ if and only if $A\intersect
p\subset q$.

With these conventions we obtain for
$x\in p\unite(r\wo q)$:
\begin{equation}\label{eq:c}
(\Delta\lk_{\set x}(s);\Delta\lk^+_{\set x}(s),\Delta\lk^-_{\set x}(s))\homeo
\begin{cases}
(\Sphere^0,\void,\void),&x\notin p,\\
(\Disk^0,\emptycmplx,\void),&x\in p\intersect q,\\
(\Disk^0,\void,\emptycmplx),&x\in p\wo r,\\
(\Disk^1,\set{1},\set{-1}),&x\in p\intersect(r\wo q).
\end{cases}
\end{equation}
We take a look only at the most interesting case.  
For $x\in p\intersect(r\wo q)$ we have
\begin{multline}\label{eq:lkxs}
\lk_{\set x}(s)=\bigl\{(\opp p,q\unite\set x,\opp r),
  (\opp{(p\wo \set x)},q\unite\set x,\opp r),
\\  (\opp{(p\wo \set x)},q,\opp r),
  (\opp{(p\wo \set x)},q,\opp{(r\wo\set x)}),
  (\opp p,q,\opp{(r\wo\set x)})
  \bigr\}
\end{multline}
with two elements comparable if and only if they are written next to
each other. Hence $\Delta\lk_{\set x}(s)$ is an interval with boundary points
$\lk_{\set x}^+(s)=\set{(\opp p,q\unite\set x,\opp r)}$ and 
$\lk_{\set x}^-(s)=\set{(\opp p,q,\opp{(r\wo\set x)})}$.

From \prettyref{eq:c} one derives inductively using \prettyref{eq:b}:
\begin{equation}\label{eq:d}
(\Delta\lk_A(s);\Delta\lk^+_A(s),\Delta\lk^-_A(s))\homeo
\begin{cases}
(\Sphere^{d_A},\void,\void),&A\intersect p=\emptyset,\\
(\Disk^{d_A},\Sphere^{d_A-1},\void),&\emptyset\ne A\intersect p\subset  q,\\
(\Disk^{d_A},\void,\Sphere^{d_A-1}),&\text{$A\intersect p\ne0$, $A\intersect p\intersect r=\emptyset$},\\
(\Disk^{d_A},\Disk^{d_A-1}_+,\Disk^{d_A-1}_-),&
    \text{$A\intersect p\intersect r\ne\emptyset$, 
          $A\intersect p\not\subset q$}
\end{cases}
\end{equation}
with $d_A=\#(A\intersect(r\wo q))+\#(A\intersect p)-1$.
Since $\lk^+(s)=\lk_{p\unite r}^+(s)$, $\lk^-(s)=\lk_{p\unite r}^-(s)$,
$\lk(s)=\lk^+(s)\unite\lk^-(s)$, the first case of~\prettyref{eq:a}
follows from the second case of~\prettyref{eq:d}, while the second case of~\prettyref{eq:a} follows from the last case of~\prettyref{eq:d}.
\end{proof}

\begin{rem}
The above Proposition will be used in the proof of \prettyref{thm:result}
to show that $\Hom_{\set{2,4}}(C_5,K_{n+2})$ is the boundary of a
regular neighbourhood of $\Delta\set{\phi\in F\Hom(C_5\wo\set{2,4},
  K_{n+2})\colon\compl\phi(3)\subset\phi(1)}$.  The proof of this
special case would require exactly the same calculation, although the
notation might possibly be more lucid.  Indeed, the poset in
\prettyref{eq:lkxs} would become
$\ind(C_5\wo\set{2,4})\wo\set\emptyset$.
\end{rem}

\begin{prop}\label{prop:nbd}
Let $P$ be the face poset of a pl-triangulation of a compact manifold $M^n$.
Let $B\subset\opp P\times P\times\opp P$ be given by
\begin{equation*}
B\deq\set{(\opp p,q,\opp r)\colon q\le r,p\intersect r\ne\emptyset,p\not\le q}.
\end{equation*}
Then $\Delta B$ is homeomorphic to the boundary of a regular neighbourhood of 
the diagonal in $M\times M$.
\end{prop}

\begin{proof}
We take up the notation of \prettyref{prop:nb-mf} and consider $B$ as
a subset of $\opp P\times\Int P$.  We also define $D\deq\set{(\opp
  p,q,\opp r)\in\opp P\times\Int P\colon p\le q}$.  $D$~is closed from
above and hence $\Delta D$ is a full subcomplex of $\Delta(\opp
P\times\Int P)$.  $(\Delta N,\Delta B)$ is the simplicial
neighbourhood of $\Delta D$.  By \prettyref{prop:nb-mf} $\Delta N$ is
a manifold with boundary $\Delta B$ and hence a regular neighbourhood
of $\Delta D$ by the Simplicial Neighbourhood
Theorem~\ref{thm:snt}.  By \prettyref{prop:int32} and
\prettyref{prop:ambient-isotopy} $\Delta D\subset\Delta(\opp
P\times\Int P)$ is ambiently isotopic to the diagonal in $M\times M$.
\end{proof}

\section{\texorpdfstring{$\Hom(C_5, K_{n+2})$}{Hom(C[5],K[n+2])}}
\label{sec:C5}
We come back to our main object of study and collect the relevant
results of the two preceding sections.  First we can replace
$\Hom(C_5, K_{n+2})$ by a smaller complex.
\begin{prop}\label{prop:C5-res-homeo}
Let $n\ge0$ and consider $\set{2,4}\subset V(C_5)$.  With the notation
of \prettyref{def:res} we have $\Hom(C_5,
K_{n+2})\homeo\Hom_{\set{2,4}}(C_5, K_{n+2})$.
\end{prop}

\begin{proof}
Since the independence complex of $C_5$ is also a $5$-gon, this
follows from \prettyref{lem:mf-homeo}.  Alternatively it is easily
checked directly that the precondition of \prettyref{lem:res-homeo} is
satisfied, and we have done so in \prettyref{ex:c5}.
\end{proof}
The graph $C_5\wo\set{2,4}$ consists of an edge and an isolated
vertex, and hence \[\Hom(C_5\wo\set{2,4},
K_{n+2})\homeo\Hom(K_2,K_{n+2})\times\Hom(K_1,K_{n+2})\homeo\Sphere^n\times\Disk^{n+1}.\]
In \prettyref{sec:Int} we have developed the techniques to describe
the subcomplex $\Hom_{\set{2,4}}(C_5, K_{n+2})$.
\begin{thm}\label{thm:result}
Let $n\ge0$.  Then $\Hom(C_5, K_{n+2})$ is homeomorphic to the Stiefel manifold
$\set{(x,y)\in\Sphere^n\times\Sphere^n\colon\kron xy=0}$.
\end{thm}

\begin{proof}
We examine
$\Hom_{\set{2,4}}(C_5,K_{n+2})$.  The face poset of this complex is
\begin{multline*}
\bigl\{\phi\colon\set{1,3,5}\to\power(V(K_{n+2}))\wo\set\emptyset
    \colon\\
    \text{$\phi(3)\intersect\phi(5)=\emptyset$,
          $\compl(\phi(1)\unite\phi(3))\ne\emptyset$,
          $\compl(\phi(5)\unite\phi(3))\ne\emptyset$}
    \bigr\}.
\end{multline*}
The first condition comes from the edge $\set{1,5}$, the second one
ensures that $\phi$ can be extended to the vertex~$2$, the last one
that it can be extended to~$4$.
If we define $P$ to be the poset 
$P\deq\power(V(K_{n+2}))\wo\set{\emptyset,V(K_{n+2})}$ then it follows that the map
\begin{equation}\label{eq:mapToPPP}
\begin{split}
F\Hom_{\set{2,4}}(C_5,K_{n+2})&\to\set{(\opp p,q,\opp r)\in\opp
  P\times P\times\opp P \colon\text{ $q\le r$, $p\not\le q$,
    $p\intersect
    r\ne\emptyset$}}\\
\phi&\mapsto(\compl\phi(3),\phi(1),\compl\phi(5))
\end{split}
\end{equation}
is an isomorphism.  Since $P$ is the face poset of the boundary of an
$(n+1)$-simplex, it follows from \prettyref{prop:nbd} that
$\Hom_{\set{2,4}}(C_5,K_{n+2})$ is homeomorphic to the boundary of a
regular neighbourhood of the diagonal in $\Sphere^n\times\Sphere^n$.
Regular neighbourhoods are unique up to isotopy and
$\set{(x,y)\in\Sphere^n\times\Sphere^n\colon\kron xy=0}$ is the
boundary of the regular neighbourhood
$\set{(x,y)\in\Sphere^n\times\Sphere^n\colon\kron xy\ge0}$ of the
diagonal $\set{(x,x)\in\Sphere^n\times\Sphere^n}$.
\end{proof}

\subsection{Involutions}
For the theorem of Babson and Kozlov mentioned in the introduction,
the $\Z_2$-action on $\Hom(C_{2r+1}, K_{n+2})$ induced by an
automorphism of $C_{2r+1}$ that flips an edge is important.  In
\cite{csorba-thesis} the conjecture that $\Hom(C_5,
K_{n+2})\homeo\set{(x,y)\in(\Sphere^n)^2\colon\kron xy=0}$ is
strengthened to a conjecture that there is a homeomorphism that takes
the $\Z_2$-action on $\Hom(C_5,K_{n+2})$ to the action given by
$(x,y)\mapsto(x,-y)$.  In \cite{codd} it was shown and used that there
is an equivariant map between these spaces.  This is essentially the
map $h$ that will appear in \prettyref{thm:equiv-res}.  We examine the
constructions presented in this work more closely to show a homotopy
version of the equivariant conjecture.  

First we strengthen \prettyref{prop:nbd}.

\begin{prop}\label{prop:reg-nbhd-diag}
In the situation of \prettyref{prop:nb-mf}, $(\Delta N, \Delta B)$ is
a regular neighbourhood of the diagonal in $M\times M$.
\end{prop}

\begin{proof}
In \prettyref{prop:nbd} we have seen that $\Delta N$ is a regular
neighbourhood of the subspace $\Delta D$ considered there and that
there is an ambient isotopy moving $\Delta D$ to the diagonal.  We
will now show that this isotopy can be chosen to be the identity
outside of~$\Delta N\wo\Delta B$.

The image of the isotopy constructed in the proof of
\prettyref{prop:isotopy} is contained in the union of the cells
$\sigma\times\sigma$ for simplices $\sigma$ of $\Delta P$.  We call
this union the fat diagonal.  Also every point in the image of the
isotopy is contained in the interior of one cell of the fat diagonal
or in $\Delta(\Int P)$.  The complex $\Delta(\opp P\times\Int P)$ of
\prettyref{prop:nb-mf} is a subdivision of $\Delta P\times\Delta P$ as
discussed in \prettyref{prop:int32}.  The subcomplex $\Delta N$ is
contained in the fat diagonal, and no simplex of $\Delta B$ meets the
interior of a cell of the fat diagonal. This proves the Proposition.
\end{proof}

\begin{thm}\label{thm:equiv-res}
Let $\phi$ be an involution on $\Hom(C_5, K_{n+2})$ induced by an
automorphism of $C_5$ that flips an edge.  Let
$V\deq\set{(x,y)\in\Sphere^n\times\Sphere^n\colon\kron xy=0}$ and
$\psi\colon V\to V$ be the involution $\phi(x,y)=(x,-y)$.  Then the
following statements hold.
\begin{enumerate}
\item
There is a homeomorphism $g\colon\Hom(C_5, K_{n+2})\xto\homeo V$ with
$\psi\cmps g\homot g\cmps\phi$.
\item
There is a homotopy equivalence $h\colon\Hom(C_5, K_{n+2})\xto\homot V$ with
$\psi\cmps h=h\cmps\phi$.
\end{enumerate}
\end{thm}

\begin{proof}
Let us consider the automorphism of $C_5$ given by $j\mapsto6-j$.
\prettyref{fig:indc5} shows that the homeomorphism
$\Hom(C_5,K_{n+2})\homeo\Hom_{\set{2,4}}(C_5, K_{n+2})$ can be chosen
to respect the induced involutions.  In \prettyref{eq:mapToPPP},
$\Hom_{\set{2,4}}(C_5, K_{n+2})$ is identified with a subspace of
$\Delta(\opp P\times\Int P)$ where $P$~is the face poset of the
boundary of the $(n+1)$-simplex.  This identification is compatible
with the involution on $\opp P\times\Int P$ given by $(\opp p,q,\opp
r)\mapsto(\opp p,a(r),\opp{a(q)})$ where $a$ is the map sending a face
to the face opposite to it.  On
$\Sphere^n\times\Sphere^n\homeo\Delta(\opp P\times\Int P)$ this map
induces the identity on the first factor and the antipodal map on the
second, and we will from now on consider the involution
$(x,y)\mapsto(x,-y)$ on $\Sphere^n\times\Sphere^n$.  By
\prettyref{prop:reg-nbhd-diag}, $\Hom_{\set{2,4}}(C_5, K_{n+2})$ is
the boundary of a regular neighbourhood of the diagonal in
$\Sphere^n\times\Sphere^n$.  The same is true of~$V$, and for both
spaces the involution on $\Sphere^n\times\Sphere^n$ exchanges the
components of their complements.  We define involutions
\begin{align*}
\Phi\colon[-1,1]\times\Hom(C_5, K_{n+2})&\to[-1,1]\times\Hom(C_5,K_{n+2}),\\
(t,x)&\mapsto(-t,\phi(x)),
\intertext{and}
\Psi\colon[-1,1]\times V&\to[-1,1]\times V,\\
(t,x)&\mapsto(-t,\psi(x)).
\end{align*}
A structure of regular neighbourhood of the diagonal with boundary $V$
defines for $\eps>0$ small enough an equivariant embedding of
$[-1,1]\times V$ in $\Sphere^n\times\Sphere^n$ which sends
$\set{-1}\times V$ to the boundary of the $\eps$-neighbourhood of the
diagonal and $\set0\times V$ to $V$.  Analogously and additionallly
using the homeomorphism $\Hom(C_5, K_{n+2})\homeo\Hom_{\set{2,4}}(C_5,
K_{n+2})$, we obtain an embedding of $[-1,1]\times\Hom(C_5, K_{n+2})$
with the same image.  This yields homeomorphisms
\begin{align*}
H\colon[-1,1]\times\Hom(C_5,K_{n+2})&\xto\homeo[-1,1]\times V,\\
g\colon\Hom(C_5, K_{n+2})&\xto\homeo V
\intertext{with}
H(-1,x)&=(-1,g(x))\text{\quad for all $x$,}\\
H\cmps\Phi&=\Psi\cmps H.
\end{align*}
It follows that
\begin{equation*}
H(1,x)=\Psi(H(-1,\phi(x)))=\Psi(-1,g(\phi(x)))=(1,\psi(g(\phi(x)))).
\end{equation*}
Therefore, if we denote by $p_V\colon[-1,1]\times V\to V$ the
projection, then $p_V\cmps H$ is a homotopy between $g$ and $\psi\cmps
g\cmps\phi$, which proves the first statement.  If we define $h(x)\deq
p_V(H(0,x))$, then $h\homot g$, since $g(x)=p_V(H(-1,x))$, and
therefore $h$ is a homotopy equivalence.  Also
\begin{align*}
h(\phi(x))&=p_V(H(0,\phi(x)))=p_V(H(\Phi(0,x)))=p_V(\Psi(H(0,x)))=\psi(p_V(H(0,x)))
\\&=\psi(h(x)),
\end{align*}
which proves the second statement.
\end{proof}

\begin{rem}
The subspaces $V$ and $\Hom_{\set{2,4}}(C_5,K_{n+2})$ of
$\Sphere^n\times\Sphere^n$ are characteristic submanifolds of the
involution on $\Sphere^n\times\Sphere^n$ in the language of
\cite{involutions}.  It might be possible to use the surgery
techniques presented there to show that the induced involutions on
them are equivalent, at least for $n\ge3$ when the $s$-cobordism
theorem is available.  On the other hand the description of
$\Hom_{\set{2,4}}(C_5,K_{n+2})$ as a subpolyhedron of
$\Sphere^n\times\Sphere^n$ that we have used is explicit enough to
hope that it might be possible to prove this by a more direct aproach.
\end{rem}

\bibliographystyle{cas-ea}
\bibliography{math,topology,combi}

\begin{thebibliography}{CLSW04}
\expandafter\ifx\csname url\endcsname\relax
  \def\url#1{{\tt #1}}\fi
\expandafter\ifx\csname urlprefix\endcsname\relax\def\urlprefix{URL }\fi
\expandafter\ifx\csname selectlanguage\endcsname\relax
  \def\selectlanguage#1{\relax}\fi

\bibitem[BK06a]{babson-kozlov-i}
{\sc Babson, E.} and {\sc Kozlov, D.~N.}
\newblock {\selectlanguage{english}Complexes of graph homomorphisms}.
\newblock {\em Isr. J. Math.\/}, {\bf 152}:285--312, 2006.
\newblock \href{http://arxiv.org/abs/math/0310056}{math.CO/0310056}.

\bibitem[BK06b]{babson-kozlov-ii}
---.
\newblock {\selectlanguage{english}Proof of the {L}ov\'asz conjecture}.
\newblock {\em Annals of Mathematics\/}, 2006.
\newblock In press, \href{http://arxiv.org/abs/math/0402395}{math.CO/0402395}.

\bibitem[CL05]{csorba-lutz}
{\sc Csorba, P.} and {\sc Lutz, F.~H.}
\newblock Graph coloring manifolds, 2005.
\newblock 22~pp., \href{http://arxiv.org/abs/math/0510177}{math.CO/0510177}.

\bibitem[CLSW04]{shore-subdiv}
{\sc Csorba, P.}, {\sc Lange, C.}, {\sc Schurr, I.}, and {\sc Wassmer, A.}
\newblock Box complexes, neighborhood complexes, and the chromatic number.
\newblock {\em Journal of Combinatorial Theory, Series A\/}, {\bf
  108}:159--168, 2004.
\newblock \href{http://arxiv.org/math/0310339}{math.CO/0310339}.

\bibitem[Cso04]{csorba-box}
{\sc Csorba, P.}
\newblock Homotopy types of box complexes, 2004.
\newblock \href{http://arxiv.org/math/0406118}{math.CO/0406118}.

\bibitem[Cso05]{csorba-thesis}
---.
\newblock {\em Non-tidy Spaces and Graph Colorings\/}.
\newblock Ph.D. thesis, ETH Z\"urich, 2005.

\bibitem[Koz05a]{kozlov-cycles}
{\sc Kozlov, D.~N.}
\newblock Cohomology of colorings of cycles, 2005.
\newblock Preprint, \href{http://arxiv.org/math/0507117}{math.AT/0507117}.

\bibitem[Koz05b]{kozlov-folds}
---.
\newblock A simple proof for folds on both sides in complexes of graph
  homomorphisms.
\newblock {\em Proceedings of the AMS\/}, 2005.
\newblock To appear, \href{http://arxiv.org/math/0408262}{math.CO/0408262}.

\bibitem[LdM71]{involutions}
{\sc L\'opez~de Medrano, S.}
\newblock {\em Involutions on Manifolds\/}, vol.~59 of {\em Ergebnisse der
  Mathematik und ihrer Grenzgebiete\/}.
\newblock Springer-Verlag, 1971.

\bibitem[Lov78]{lovasz}
{\sc Lov\'asz, L.}
\newblock Kneser's conjecture, chromatic number and homotopy.
\newblock {\em J. Combinatorial Theory, Ser. A\/}, {\bf 25}:319--324, 1978.

\bibitem[MZ]{matousek-ziegler}
{\sc Matou{\v s}ek, J.} and {\sc Ziegler, G.~M.}
\newblock Topological lower bounds for the chromatic number: {A} hierarchy.
\newblock {\em Jahresbericht der DMV\/}.
\newblock To appear, \href{http://arxiv.org/math/0208072}{math.CO/0208072}.

\bibitem[RS82]{rourke-sanderson}
{\sc Rourke, C.~P.} and {\sc Sanderson, B.~J.}
\newblock {\em Introduction to Piecewise-Linear Topology\/}.
\newblock Springer-Verlag, 1982.

\bibitem[Sch05]{codd}
{\sc Schultz, C.}
\newblock {\selectlanguage{english}A short proof of
  {$w_1^n(\Hom(C_{2r+1},K_{n+2}))=0$} for all~$n$ and a graph colouring theorem
  by {B}abson and {K}ozlov}, 2005.
\newblock Preprint, 8pp.,
  \href{http://arxiv.org/abs/math/0507346}{math.AT/0507346}.

\bibitem[Sch06]{hom-loop}
---.
\newblock Graph colourings, spaces of edges and spaces of circuits, 2006.
\newblock Preprint, \href{http://arxiv.org/math/0606763}{math.CO/0606763}.

\bibitem[Seg68]{segal-classifying-spaces}
{\sc Segal, G.}
\newblock Classifying spaces and spectral sequences.
\newblock {\em Publ. Math. IHES\/}, {\bf 34}:105--112, 1968.

\bibitem[Wal88]{walker-posets}
{\sc Walker, J.~W.}
\newblock Canonical homeomorphisms of posets.
\newblock {\em European J. Combin.\/}, {\bf 9}(2):97--107, 1988.
\newblock ISSN 0195-6698.

\bibitem[Zee63]{unknot}
{\sc Zeeman, E.}
\newblock Unknotting of combinatorial balls.
\newblock {\em Annals of Mathematics\/}, {\bf 78}:501--526, 1963.

\bibitem[{\v{Z}}iv05]{wiposets}
{\sc {\v{Z}}ivaljevi\'c, R.~T.}
\newblock {$WI$}-posets, graph complexes and {$\Z_2$}-equivalences.
\newblock {\em Journal of Combinatorial Theory, Series A\/}, {\bf
  111}:204--223, 2005.
\newblock \href{http://www.arxiv.org/math/0405419}{math.CO/0405419}.

\end{thebibliography}

\end{document}